\documentclass{article}
\usepackage{amssymb}
\usepackage{amsmath}
\usepackage{latexsym}
\usepackage{mathrsfs}
\usepackage{graphicx}
\usepackage[all]{xy}
\newtheorem{Th}{Theorem}
\newtheorem{Lemma}[Th]{Lemma}

\def\into{\mathop{\hookrightarrow}\limits}
\def\onto{\mathop{\twoheadrightarrow}\limits}
\title{A short proof of Spectral Sequences of Filtered Complexes}
\author{XIONG Rui}
\begin{document}

\def\E{\mathsf{E}}
\def\B{\mathsf{B}}
\def\Z{\mathsf{Z}}
\def\H{\mathsf{H}}
\def\F{\mathcal{F}}
\def\L{\mathcal{L}}
\def\di{\textsf{d}}
\def\im{\operatorname{im}}
\def\fref#1{}

\maketitle

\begin{abstract}
The purpose of this article is to present my new proof of the the construction and the convergence theorem of spectral sequences of filtered complexes, which is much shorter and cleaner than the ``standard'' proof.
For the historical remarks, see the last section.
\end{abstract}

Thanks to Yuri V. Volkov for finding a mistake in the original proof.

\section{Some fundamental Lemmate}

Here I present several more or less trivial lemmate without proofs which I will use in the proof.
This are all common knowledge in algebra, see for example \cite{zbMATH01703931}.

\begin{Lemma}[Modular property]\label{ModPro}
Let $A,B,C$ be three subgroups of some bigger abelian group, if $A\subseteq C$, then
$$(A+B)\cap C=A+(B\cap C).$$
As a result, it makes no doubt to write $A+B\cap C$.
\end{Lemma}

\begin{Lemma}[Exchange Limit]\label{ExLim}
Let $C_\bullet$ be a directed family of submodules of some bigger abelian group, that is,
each pair of $C_{i},C_j$ are submodules of some $C_k$. Assume $A\subseteq B$, then we have
$$\begin{array}{c}
\bigcup  \big(A+C_\bullet \cap B\big)=A+\big(\bigcup  C_\bullet \big) \cap B.
\end{array}$$
If furthermore, $C_\bullet$ is  bounded below, that is, some $C_i$ equals to $\bigcap C_\bullet$, then,
$$\begin{array}{c}
\bigcap \big(A+C_\bullet \cap B\big)=A+\big(\bigcap C_\bullet \big) \cap B.
\end{array}$$
\end{Lemma}

\begin{Lemma}[Zassenhaus' Butterfly Lemma]\label{Btf}
For four subgroups $A,B,C,D$ of some bigger abelian group, if $A\subseteq B$ and $C\subseteq D$, then
$$\frac{A+D\cap B}{A+C\cap B}\cong \frac{C+B\cap D}{C+A\cap D}.$$
\end{Lemma}


\section{Functorial Zassenhaus' Butterfly Lemma}

Let $X\stackrel{f}\to X'$ be a homomorphism between abelian groups, let 
\begin{itemize}
\item $A,B,C,D$ be four subgroups of $X$ with $A\subseteq B$ and $C\subseteq D$, 
\item $A',B',C',D'$ be four subgroups of $X'$ with $A'\subseteq B'$ and $C'\subseteq D'$, and
\item $f(A)\subseteq A', f(B)\subseteq B', f(C)\subseteq C', f(D)\subseteq D'$. 
\end{itemize}

\begin{Th}[Functorial Butterfly] Then $f$ induces
$$\hat{f}: \frac{A+D\cap B}{A+C\cap B}\to \frac{A'+D'\cap B'}{A'+C'\cap B'}. $$
\end{Th}

\begin{Lemma}\label{surjlemma}If $f(B\cap D)=B'\cap D'$, then $\hat{f}$ is surjective.
\end{Lemma}
Since any element on the right is presented by an element $x$ of $D'\cap B'$ (since $A'$ is in the denominator).
By the assumption, $f(y)=x$ for some $y\in D$, now $y\in f^{-1}(B')=B$, so $x$ is mapped by $y$.

%

\begin{Lemma}\label{injlemma}If $A+C=f^{-1}(A'+C')$, then $\hat{f}$ is injective. 
\end{Lemma}
Pick any element on the left, assume it is presented by $x\in D\cap B$, 
assume $f(x)\in (A'+C')\cap B'$, now $x\in A+C$ by assumption, 
so is in the denominator $(A+C)\cap B$. 
%
%

The above two lemmate also follow immediately from theorem \ref{ExtBuf} as a remark. 

\begin{Th}[Adjointness]\label{adjointofButterfly} If $A=f^{-1}(A')$, $B=f^{-1}(B')$ and $f(C)=C'$, $f(D)=D'$, then $\hat{f}$ is an isomorphism. 
\end{Th}
Firstly, $f(B\cap D)=B'\cap D'$. $x\in B'\cap D'\subseteq D'=f(D)$, 
there exists $y\in D$ such that $f(y)=x\in B'$, so $y\in B$ also as desired. 
Secondly, $A+C=f^{-1}(A'+C')$. If $f(x)=y+f(z)$ with $y\in A'$ and $z\in C$ then $f(x-z)=y\in A'$, i.e. $x-z\in A$. 
Therefore, $x\in A+C$ as desired again. 


\section{The Construction}


Let $(C_\bullet,\di)$ be a complex with filtration $\F$ which is exhaustive (i.e. $C_n=\bigcup_{p\in \mathbb{Z}} \F_pC_n$) and
bounded below (i.e. for each $n$, $p\ll 0\Rightarrow\F_pC_n=0$). Thus the filtration automatically satisfy the condition in lemma \ref{ExLim}.
For convenience, we also denote
$$\F_{-\infty}C_n= \bigcap_{p}F_pC_{n}=0,\qquad \F_{\infty}C_n=\bigcup_{p}F_pC_n=C_n.$$
To invest its homology, define
$$\F_pK_{n}=\di^{-1}(\F_pC_{n-1})\quad\textrm{ and }\quad \F_pI_{n}=\di(\F_pC_{n+1})$$
for each $p\in \mathbb{Z}\sqcup\{\infty\}$ and $n$.
We now use $\{\F_pK_n\}_p$ and $\{\F_pI_n\}_p$ to refine the filtration $\{\F_pC_n\}_p$. Define
$$\begin{cases}
\Z_{pq}^r=\F_{p-1}C_{p+q}+\F_{p-r}K_{p+q}\cap \F_{p}C_{p+q},\\
\B_{pq}^r=\F_{p-1}C_{p+q}+\F_{p+r-1}I_{p+q}\cap \F_{p}C_{p+q},
\end{cases}\textrm{ with }r=0,1,\ldots, +\infty.$$
See lemma \ref{ModPro} for the meaning for notations.
Now, define
$$\E_{pq}^r=\Z_{pq}^r/\B_{pq}^r.$$
We have
$$\begin{array}{rll}
\dfrac{\Z_{pq}^r}{\Z_{pq}^{r+1}}&=
\dfrac{\F_{p-1}C_{p+q}+F_{p-r}K_{p+q}\cap \F_{p}C_{p+q}}
{\F_{p-1}C_{p+q}+\F_{p-r-1}K_{p+q}\cap \F_{p}C_{p+q}}& \because \textrm{definition}\\ [2ex]
&=\displaystyle
\dfrac{\F_{p-r-1}K_{p+q}+ \F_{p}C_{p+q}\cap \F_{p-r}K_{p+q}}
{\F_{p-r-1}K_{p+q}+ \F_{p-1}C_{p+q}\cap \F_{p-r}K_{p+q}} & \because \textrm{butterfly lemma } \ref{Btf}\\ [2ex]
&\displaystyle\stackrel{\di}\to
\dfrac{\F_{p-r-1}C_{p+q-1}+ \F_{p}I_{p+q-1}\cap \F_{p-r}C_{p+q-1}}
{\F_{p-r-1}C_{p+q-1}+ \F_{p-1}I_{p+q-1}\cap \F_{p-r}C_{p+q-1}}& \because \textrm{theorem }\ref{adjointofButterfly}
\\ [2ex]
& =\dfrac{\B_{p-r,q+r-1}^{r+1}}{\B_{p-r,q+r-1}^{r}} &\because \textrm{definition}.
\end{array}$$
Then we can define
$$\di^r_{pq}=\left[\E_{pq}^r=\frac{\Z_{pq}^r}{\B_{pq}^r}\onto \frac{\Z_{pq}^r}{\Z_{pq}^{r+1}}\cong
\frac{\B_{p-r,q+r-1}^{r+1}}{\B_{p-r,q+r-1}^{r}}\into \frac{\Z_{p-r,q+r-1}^{r}}{\B_{p-r,q+r-1}^{r}}
=\E_{p-r,q+r-1}^r\right]$$
clearly, $\di^r\circ \di^r=0$, and at the position $(p,q)$ of stage $r$,
$$\ker \di^r = \frac{\Z_{pq}^{r+1}}{\B_{pq}^r},\quad \im \di^r =\frac{\B_{pq}^{r+1}}{\B_{pq}^r}, \quad\textrm{and} \quad \frac{\ker\di^r }{\im \di^r}=\E_{pq}^{r+1}. $$
Thus the above construction of $\{E_{pq}^r\}$ give rise to a spectral sequence.

\section{Convergence Theorem}
\begin{Th}[Convergence Theorem]\label{ConvergenceTh}
Given a complex $C_\bullet$ with filtration $\F$ which is exhaustive and bounded below,
then the $\{\E_{pq}^r\}$ constructed above is convergent to $\H_n(C_\bullet)$.
\end{Th}
\textsc{Proof. }
Using butterfly lemma  again,
$$\begin{array}{rll}
\E_{pq}^\infty = \dfrac{\bigcap_r \Z_{pq}^r}{\bigcup_r \B_{pq}^r}
& =\dfrac{\bigcap_{r}(\F_{p-1}C_{p+q}+\F_{p-r}K_{p+q}\cap \F_{p}C_{p+q})}
{\bigcup_{r} (\F_{p-1}C_{p+q}+\F_{p+r-1}I_{p+q}\cap \F_{p}C_{p+q})}&\because \textrm{definition}\\ [2ex]
&=\dfrac{\F_{p-1}C_{p+q}+\F_{-\infty}K_{p+q}\cap \F_{p}C_{p+q}}
{\F_{p-1}C_{p+q}+\F_{+\infty}I_{p+q}\cap \F_{p}C_{p+q}}& \because\textrm{lemma } \ref{ExLim}\\ [2ex]
&=\dfrac{\F_{p-1}C_{p+q}+\ker \di \cap \F_{p}C_{p+q}}
{\F_{p-1}C_{p+q}+\im \di \cap \F_{p}C_{p+q}} &\because \textrm{definition} \\ [2ex]
&=\dfrac{\im\di+\F_{p}C_{p+q}\cap \ker \di}{\im\di+\F_{p-1}C_{p+q}\cap \ker \di} &\because \textrm{butterfly lemma } \ref{Btf}.
\end{array}$$
Denote $\F_p\H_{p+q}=\im\di+\F_{p}C_{p+q}\cap \ker \di$, then
$$\begin{array}{l}
\bigcup_{p}\F_p\H_n  = \im\di+\big(\bigcup_{p}\F_{p}C_{n}\big)\cap \ker \di = \im\di +C_n\cap \ker\di = \ker \di, \\
\bigcap_{p}\F_p\H_n  =\im\di+\big(\bigcap_{p}\F_{p}C_{n}\big)\cap \ker \di  = \im\di +0\cap \ker\di = \im\di. 
\end{array}$$
from lemma \ref{ExLim} and the assumptions. 
So, we define a filtration over $\H_n(C_\bullet)=\dfrac{\ker \di}{\im \di}$ whose subquotient coincide with $\E_{pq}^\infty$.
The proof is complete.

\section{Computation of $E_{pq}^0$}

Firstly,
$$\begin{array}{rll}
\Z_{pq}^{0}& =\F_{p-1}C_{p+q}+\F_{p}K_{p+q}\cap \F_{p}C_{p+q}\\
&= \F_{p-1}C_{p+q}+\di^{-1}(\F_{p}C_{p+q-1})\cap \F_{p}C_{p+q}\\
& =\F_{p-1}C_{p+q},
\end{array}$$
since $\F_{p}C_{p+q}\subseteq\di^{-1}(\F_{p}C_{p+q-1})$; and
$$\begin{array}{rll}
\B_{pq}^{0} & =\F_{p-1}C_{p+q}+\F_{p-1}I_{p+q}\cap \F_{p}C_{p+q} \\
& =\F_{p-1}C_{p+q}+\di(\F_{p-1}C_{p+q+1})\cap \F_{p}C_{p+q}\\
& = \F_{p}C_{p+q},
\end{array}$$
since $\di(\F_{p-1}C_{p+q+1})\subseteq \F_{p-1}C_{p+q}$.
Hence $E_{pq}^0=\dfrac{\Z_{pq}^{0}}{\B_{pq}^{0}}=\dfrac{\F_{p}C_{p+q}}{\F_{p-1}C_{p+q}}$.
Note that the isomorphism $\dfrac{\Z_{pq}^0}{\Z_{pq}^1}\to \dfrac{\B_{p,q-1}^1}{\B_{p,q-1}^0}$ is induced by $\di$,
and by definition $\di^0$ factors through the isomorphism, so it is induced by $\di$.


\section{Remarks on Zassenhaus' Butterfly lemma}
If we reflect the ``standard'' proof in the textbooks, the reason making the it painful is 
that what we did in the proof is exactly repeating the proof of Zassenhaus' butterfly lemma element by element.

%

The reason that Zassenhaus' lemma is called butterfly is because in the proof we use the following ``butterfly''.
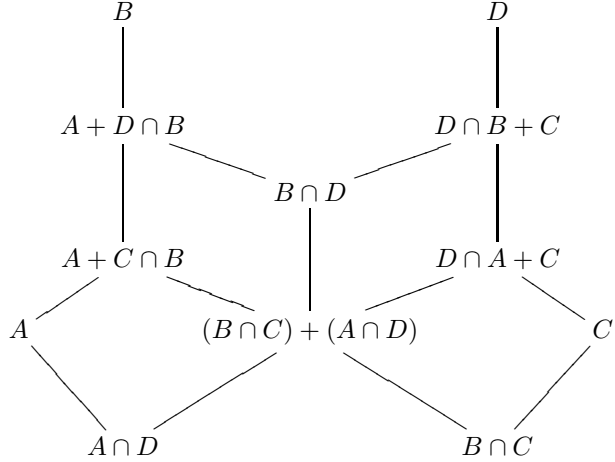
$$\xymatrix@C=0pc@R=1pc{
&&B\ar@{-}[dd]&&D\ar@{-}[dd]\\ \\
&&A+D\cap B\ar@{-}[dr]\ar@{-}[dd]  &&D\cap B+ C\ar@{-}[dl]\ar@{-}[dd]\\
&&& B\cap D\ar@{-}[dd] & \\
&& A+C\cap B\ar@{-}[lld]\ar@{-}[rd] &&D\cap A+ C\ar@{-}[rrd]\ar@{-}[ld]\\
A\ar@{-}[ddrr]&&&(B\cap C)+(A\cap D)\ar@{-}[ddr]\ar@{-}[ddl]&& &C\ar@{-}[lldd]\\ \\
&&A\cap D&&B\cap C&
}$$
Without the middle term, the element by element proof will not be easy, especially in our case where we have too much modules.

But actually, we have dual butterfly also. 

\begin{Th}[Extended Zassenhaus' lemma]\label{ExtBuf}
 For four subgroups $A,B,C,D$ of some bigger abelian group, if $A\subseteq B$ and $C\subseteq D$, 
then we have the following isomorphisms
$$\dfrac{A+D\cap B}{A+C\cap B}\cong  \dfrac{B\cap D}{(B\cap C)+(A\cap D)} \cong \dfrac{(A+D)\cap (B+C)}{A+C}\cong \dfrac{C+B\cap D}{C+A\cap D}$$
\end{Th}
Note that this implies lemma \ref{surjlemma} and \ref{injlemma} immediately. 

The one more isomorphism is absolutely the same to the standard proof. The diagram is as following 
$$
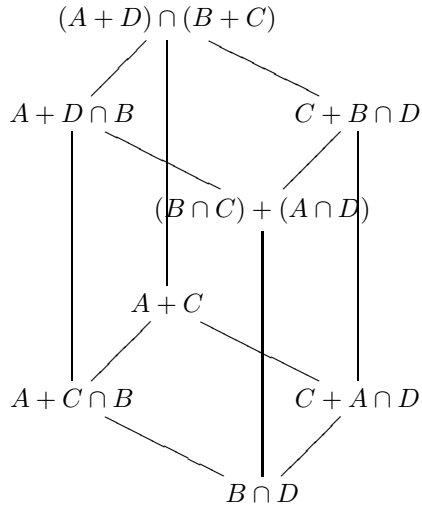
\xymatrix@!=1pc{
& (A+D)\cap (B+C)\ar@{-}[ddd]\ar@{-}[dl]\ar@{-}[drr] \\
A+D\cap B\ar@{-}[ddd]\ar@{-}[drr] &&& C+B\cap D \ar@{-}[ddd]\ar@{-}[dl]\\ 
&& (B\cap C)+(A\cap D)\ar@{-}[ddd]\\
&A+C\ar@{-}[dl]\ar@{-}[drr]& \\
A+C\cap B\ar@{-}[drr] &&& C+A\cap D\ar@{-}[dl] \\
&& B\cap D
}$$

\section{Historical Remarks}

For the readers never read the ``standard'' proof in the textbooks,
the presentation here is strange and confusing. Hence some historical remarks should be given.
As introduced in \cite{zbMATH01565334},
\begin{quote}
  Theorem 2.6 (i.e. theorem \ref{ConvergenceTh} here) first appeared in the work of [Koszul47] and [Cartan48]
who had extracted the algebraic essence underlying the work of [Leray46] on
sheaves, homogeneous spaces, and fibre spaces.
\end{quote}

In \cite{zbMATH01565334} page 34 the proof of theorem 2.6, it claims 4 pages, with huge calculation of identities between submodules, and dazzling indices.
In  \cite{zbMATH00841435} page 133 construction 5,4,6 and lemma 5.4.7, leaving the main process an exercise.
In \cite{zbMATH05306055} page 626 theorem 10.14, the author used the exact couple to avoid direct proof.
In the famous text book about spectral sequences \cite{zbMATH01565334}, the author wrote
\begin{quote}
Though the guts of our black box are laid bare, such close
examination may reveal details too fine to be enlightening. Thus we place the
proof of the theorem in a separate section and recommend skipping it to the
novice and to the weak of interest.
\end{quote}

When I was undergraduate in 2018,
I tried to understand the ``standard'' proof of this fact but failed. So I spent approximately one week to think up a ``natural'' proof, especially without element picking.
Here is what I got. The main process is by one punch, and even though all the details are written down, it is still in 2 pages leaving much margin.
Through my proof, we can also see that spectral sequence for filtered complex is in some way a kind of ``calculating homological group'' in two different ways.
Furthermore, the assumption that the complex is exhaustive and bounded is just used to exchange the order as in lemma \ref{ExLim}.

\bibliographystyle{plain}
\bibliography{bibfile}

\begin{thebibliography}{1}

\bibitem{zbMATH01703931}
Serge {Lang}.
\newblock {\em {Algebra. 3rd revised ed.}}, volume 211.
\newblock New York, NY: Springer, 3rd revised ed. edition, 2002.

\bibitem{zbMATH01565334}
John {McCleary}.
\newblock {\em {A user's guide to spectral sequences. 2nd ed.}}, volume~58.
\newblock Cambridge: Cambridge University Press, 2nd ed. edition, 2001.

\bibitem{zbMATH05306055}
Joseph~J. {Rotman}.
\newblock {\em {An introduction to homological algebra. 2nd ed.}}
\newblock Berlin: Springer, 2nd ed. edition, 2009.

\bibitem{zbMATH00841435}
Charles~A. {Weibel}.
\newblock {\em {An introduction to homological algebra. 1st pbk-ed. 1st
  pbk-ed.}}, volume~38.
\newblock Cambridge: Cambridge Univ. Press, 1st pbk-ed. edition, 1995.

\end{thebibliography}

\end{document}